\documentclass[a4paper,oneside,10pt]{article} 

\title{AN EXPONENTIAL INEQUALITY FOR SYMMETRIC RANDOM VARIABLES}
\author{Rapha\"el Cerf & Matthias Gorny}
\date{\textit{MAJ: October 15th, 2014}}

\usepackage[english]{babel}
\usepackage{amsfonts}
\usepackage{cite}
\usepackage{amsthm}
\usepackage{amssymb}
\usepackage{amsmath}
\usepackage{dsfont}

\newcommand{\R}{\mathbb{R}}
\renewcommand{\P}{\mathbb{P}}
\newcommand{\E}{\mathbb{E}}
\newtheorem*{theo}{Theorem}

\begin{document}

\renewcommand{\contentsname}{Contents}
\renewcommand{\refname}{\textbf{References}}
\renewcommand{\abstractname}{Abstract}

\begin{center}
\begin{LARGE}
AN EXPONENTIAL INEQUALITY\smallskip 

FOR SYMMETRIC RANDOM VARIABLES
\end{LARGE}\bigskip \bigskip \bigskip

\begin{large}
{\setlength{\tabcolsep}{15pt}
\begin{tabular}{p{3.4cm}p{1cm}p{3.8cm}}
\centering Rapha\"el Cerf & \centering and &\centering Matthias Gorny\tabularnewline
\centering {\it ENS Paris} 
&& 
\centering{\it Universit\'e Paris Sud}\tabularnewline
\centering && \centering{\it and ENS Paris}
\end{tabular}}

\end{large}
\end{center}

%

\bigskip\bigskip


A central question in the study of sums of independent identically distributed random variables is to control the probabilities that they deviate from their typical values. We present here a very simple exponential inequality. This inequality holds for any symmetric distribution, and does not require any integrability condition.

\begin{theo} Let $n\geq 1$ and let $X_1,\dots,X_n$ be $n$ independent identically distributed symmetric real-valued random variables.
For any $x,y>0$, we have
\[\P\big({X_1+\dots+X_n}\geq x,\,
{X_1^2+\dots+X_n^2}\leq y\big)< \exp\left(-\frac{x^2}{2y}\right)\,.\]
\end{theo}

\noindent If we apply this inequality with $nx,ny$ instead of $x,y$, we obtain an inequality controlling the first two empirical moments of the sequence $X_1,\cdots,X_n$:
\[\P\left(\frac{X_1+\dots+X_n}{n}\geq x,\,\frac{X_1^2+\dots+X_n^2}{n}\leq y\right)< \exp\left(-\frac{nx^2}{2y}\right)\,.\]
With the classical theory of large deviations \cite{DZ}, we usually obtain exponential inequalities of this type, but unfortunately they are valid for $n$ large and it is a difficult task to quantify how large $n$ has to be. During the last decades, probabilists and statisticians have been trying to find non asymptotic inequalities, which are valid for any $n\geq 1$ and with explicit constants \cite{ConcentrationBLM}. Most of these inequalities are based on the phenomenon of concentration of measure and they require a strong control on the tail of the distribution of $X_1$, typically the existence of an exponential moment. The inequality we present here is a deviation inequality which is valid for all $n\geq 1$. We suppose that the distribution is symmetric, and no further integrability condition is required. Initially, this inequality was obtained through classical results of large deviations within Cram\'er's theory and an inequality on the rate function derived to study a Curie--Weiss model of self--organized criticality \cite{CerfGorny}. Since the statement of the exponential inequality is very simple, we looked for an elementary proof, and it is this proof we present here.\medskip

\noindent Let us prove the theorem. We suppose that $\P(X_1=0)<1$, otherwise the inequality of the theorem is immediate. Let $n\geq 1$ and $x,y>0$. We set
\[S_n=X_1+\dots+X_n\,,\qquad T_n=X_1^2+\dots+X_n^2\,.\]
Let $s,t>0$.
We have
\begin{align*}
\P\big({S_n}\geq x,\,{T_n}\leq y \big)&=\P\left(sS_n\geq sx
,\,-tT_n\geq -ty\right)\\
&\leq\P\big(sS_n-tT_n\geq sx-ty\big)\\
&\leq\P\Big(\exp\big(sS_n-tT_n\big)\geq \exp\big(sx-ty\big)\Big)\,.
\end{align*}
We recall one of the most classical stochastic inequality.\medskip

\noindent{\bf Markov's inequality.} If $X$ is a non--negative random variable, then
\[\forall \lambda>0 \qquad \P(X\geq \lambda)\leq\frac{\E(X)}{\lambda}\,.\]\smallskip

\noindent Using Markov's inequality and the fact that $X_1,\dots,X_n$ are i.i.d., we get
\begin{align*}
\P\big({S_n}\geq x,\,{T_n}\leq y \big)&\leq\,
\exp\big(-sx+ty\big)\,\E\Bigg(\prod_{i=1}^n\exp\big(sX_i-tX_i^2\big)\Bigg)\\
&
=\,\exp\big(-sx+ty\big)\,\Bigg(\E\Big(\exp\big(sX_1-tX_1^2\big)\Big)
\Bigg)^n\,.
\end{align*}
The distribution of $X_1$ is symmetric thus
\begin{multline*}
\E\Big(\exp\big(sX_1-tX_1^2\big)\Big)\,=\,
\E\Big(\exp\big(-sX_1-tX_1^2\big)\Big)\\
\,=\,\frac{1}{2}\Big(
\E\Big(\exp\big(sX_1-tX_1^2\big)\Big)+
\E\Big(\exp\big(-sX_1-tX_1^2\big)\Big)\Big)\\
\,=\,\E\Big(\mathrm{cosh}(sX_1)\,\exp\left(-tX_1^2\right)\Big)\,.
\end{multline*}
We choose now $t=s^2/2$. 
We have the inequality
$$
\forall u\in \R\backslash\{0\} \qquad
\mathrm{cosh}(u)\,\exp\left(-u^2/2\right)< 1\,.$$
Since $\P(X_1=0)<1$, the above inequality implies that
$$\E\Big(\mathrm{cosh}(sX_1)\,\exp\left(-s^2X_1^2/2\right)\Big)\,<\,1\,,$$
whence also
\begin{align*}
\P\big({S_n}\geq x,\,{T_n}\leq y \big)\,<\,
\exp\big(-sx+s^2y/2\big)\,.
\end{align*}
We finally choose $s=x/y$ and we obtain the desired inequality.

\bibliographystyle{plain}
\bibliography{biblio}

\begin{thebibliography}{1}

\bibitem{ConcentrationBLM}
St{\'e}phane Boucheron, G{\`a}bor Lugosi, and Pascal Massart.
\newblock {\em {Concentration Inequalities: A Nonasymptotic Theory of
  Independence}}.
\newblock Oxford University Press, 2013.

\bibitem{CerfGorny}
Rapha{\"e}l Cerf and Matthias Gorny.
\newblock {A {C}urie-{W}eiss model of {S}elf-{O}rganized {C}riticality}.
\newblock {\em The Annals of Probability}, to appear, 2013.

\bibitem{DZ}
Amir Dembo and Ofer Zeitouni.
\newblock {\em {Large deviations techniques and applications}}, volume~38 of
  {\em {Stochastic Modelling and Applied Probability}}.
\newblock Springer-Verlag, 2010.

\end{thebibliography}

\end{document}